# A Two-Stage Approach for Combined Heat and Power Economic Emission Dispatch: Combining Multi-Objective Optimization with Integrated Decision Making

Yang Li [a,b*], Jinlong Wang [a], Dongbo Zhao [b], Guoqing Li [a], Chen Chen [b]

[a] School of Electrical Engineering, Northeast Electric Power University, Jilin 132012, China

[b] Energy Systems Division, Argonne National Laboratory, Lemont, IL 60439, USA

**ABSTRACT:** To address the problem of combined heat and power economic emission dispatch (CHPEED), a two-stage approach is proposed by combining multi-objective optimization (MOO) with integrated decision making (IDM). First, a practical CHPEED model is built by taking into account power transmission losses and the valve-point loading effects. To solve this model, a two-stage methodology is thereafter proposed. The first stage of this approach relies on the use of a powerful multi-objective evolutionary algorithm, called $\theta$-dominance based evolutionary algorithm ($\theta$-DEA), to find multiple Pareto-optimal solutions of the model. Through fuzzy c-means (FCM) clustering, the second stage separates the obtained Pareto-optimal solutions into different clusters and thereupon identifies the best compromise solutions (BCSs) by assessing the relative projections of the solutions belonging to the same cluster using grey relation projection (GRP). The novelty of this work is in the incorporation of an IDM technique FCM-GRP into CHPEED to automatically determine the BCSs that represent decision makers' different, even conflicting, preferences. The simulation results on three test cases with varied complexity levels verify the effectiveness and superiority of the proposed approach.

**Keywords:** Cogeneration; economic emission dispatch; two-stage approach;

---

[*]Corresponding author. E-mail address: liyang@neepu.edu.cn (Yang Li).

multi-objective optimization; integrated decision making; valve-point loading effects; $\theta$-dominance based evolutionary algorithm; grey relational projection; integrated energy system.

**NOMENCLATURE**

**Acronyms**

| | |
|---|---|
| CHP | combined heat and power |
| CHPED | CHP economic dispatch |
| CHPEED | CHP economic emission dispatch |
| MOEAs | multi-objective evolutionary algorithms |
| MOPSO | multi-objective particle swarm optimization |
| EFA | enhanced firefly algorithm |
| MOO | multi-objective optimization |
| NSGA-II | non-dominated sorting genetic algorithm-II |
| MLCA | multi-objective line-up competition algorithm |
| TVA-PSO | time-varying acceleration particle swarm optimization |
| NLP | non-linear programming |
| NBI | normal boundary intersection |
| POF | Pareto-optimal front |
| BCSs | best comprise solutions |
| VPLE | valve-point loading effects |
| $\theta$-DEA | $\theta$-dominance based evolutionary algorithm |
| IDM | integrated decision making |
| FCM | fuzzy c-means |
| GRP | grey relation projection |
| FOR | feasible operation region |
| SPEA 2 | strength Pareto evolutionary algorithm 2 |
| RCGA | real coded genetic algorithm |
| PSO | particle swarm optimization |
| EP | evolutionary programming |
| AIS | artificial-immune system optimization |
| DE | differential evolution |
| BCO | bee colony optimization |
| SO | scalar optimization |

**Symbols**

| | |
|---|---|
| $P_i^{\min}$ | lower limit of the output of unit $i$ |
| $e_i, \zeta_i$ | cost change coefficients caused by the VPLE of power-only unit $i$ |
| $C_{Total}$ | total fuel cost |
| $C_{p,i}$ | fuel cost of power-only unit $i$ |
| $C_{c,j}$ | fuel cost of CHP unit $j$ |
| $C_{h,k}$ | fuel cost of heat-only unit $k$ |
| $a_i, b_i, d_i$ | cost coefficients of power-only unit $i$ |
| $\varphi_k, \eta_k, \lambda_k$ | cost coefficients of heat-only unit $k$ |
| $\alpha_j, \beta_j, \gamma_j, \delta_j, \varepsilon_j, \xi_j$ | cost coefficients of CHP unit $j$ |
| $P_i, O_j$ | generating capacities of power-only unit $i$ and CHP unit $j$ |
| $H_j, T_k$ | calorific values of the CHP unit $j$ and heat-only unit $k$ |
| $N_p, N_h, N_c$ | numbers of power-only units, heat-only units and CHP units |
| $E_{Total}$ | total emission of polluting gases |
| $E_S$ | total emission of $SO_2$ and $NO_x$ |
| $E_C$ | emission of $CO_2$ |
| $\mu_i, \kappa_i, \pi_i, \sigma_i, \upsilon_i$ | $SO_2$ and $NO_x$ emission coefficients of power-only unit $i$ |
| $\tau_j$ | $SO_2$ and $NO_x$ emission coefficients of CHP unit $j$ |
| $\rho_k$ | $SO_2$ and $NO_x$ emission coefficients of heat-only unit $k$ |
| $\vartheta_i, \psi_j, \varpi_k$ | $CO_2$ emission coefficients of different types of units |
| $P_D, P_L$ | power demand and power transmission loss |
| $B_{i,j}$ | loss coefficient related to the productions of unit $i$ and $j$ |
| $B_{0,i}$ | loss coefficient concerned with the production of unit $i$ |
| $B_{0,0}$ | loss coefficient parameter |
| $H_D$ | total system heat demand |
| $t$ | a time interval |
| $N_T$ | number of time intervals |
| $UR_i, DR_i$ | ramp-up and ramp-down rate limits of unit $i$ |
| $\tilde{F}(x)$ | normalized objective vector of solution $x$ |
| $\Lambda$ | set of reference points |
| $N_d$ | number of reference points |
| $\lambda_j$ | $j$th reference point |
| $L$ | a line passing through the origin and reference point $\lambda_j$ |
| **u** | projection of $\tilde{F}(x)$ on **L** |
| $\theta$ | penalty parameter |
| $N_{cl}$ | number of clusters in the population |
| $Dis_{j,1}(x)$ | distance between the origin and **u** |
| $Dis_{j,2}(x)$ | perpendicular distance between $\tilde{F}(x)$ and **L** |
| $A_0$ | initial population |
| $Z^*$ | ideal initial point |
| $Z^{nad}$ | worst point |
| $Iter$ | current iteration number |
| $A$ | set of nondominated solutions |
| $Q_{Iter}, R_{Iter}, A_{Iter}$ | offspring, new and current populations |
| $F_i$ | $i$th Pareto non-dominated level |
| $J$ | loss function |
| $J_{cur}, J_{pre}$ | current and previous loss function $J$ |
| $U$ | membership degree matrix |
| $V$ | cluster centers |
| $W$ | input vector |
| $\mu_{ij}$ | membership degree |
| $m$ | parameter that determines the degree of fuzziness |
| $Prj_l$ | projection of dispatch scheme $l$ onto the ideal scheme |
| $N_g$ | number of indicators used to assess a scheme |
| $Grc$ | grey relation coefficient |
| $r_i$ | weight of each indicator in a scheme |
| $\Gamma$ | set of targeted points |
| $Edis_1, Edis_2$ | Euclidean distances between extreme and boundary solutions |
| $\overline{Edis}$ | average value of all distances |
| $\varepsilon$ | preset threshold |
| "+"/"−" | positive/negative schemes |

## 1. Introduction

With the increasing energy crisis and environmental issues, combined heat and power (CHP) generation, also called cogeneration, has attracted ever-growing concerns in recent years and it has also proven to be an effective way for addressing these challenges [1]. In traditional thermal power plants, a lot of thermal energy is wasted without conversion into electricity during power generation. Even in terms of the most advanced combined cycle power plant, the energy conversion efficiency is by far only in the range from 50% to 60% [2]. The central and most fundamental principle of cogeneration is to improve the total energy conversion efficiency by recovering and reutilizing the waste heats in the energy conversion process [3], and thereby the fuel utilization efficiency of CHP units can achieve 90% and above [4]. At the same time, compared with traditional power-only units and heat-only units, CHP units can save 10%~40% of the cost of generation, which means that less fuels are needed to produce equal amounts of heat and electricity [5]. Furthermore, recent research suggest that CHP units are considered as an environmentally friendly system, since the greenhouse gas emissions can be reduced by nearly 13%~18% by making use of cogenerations [6].

CHP economic dispatch (CHPED) has been recognized as an important means to achieve optimal operation for CHP systems, since it is able to significantly reduce the unit energy consumption of coal-fired power plants through optimizing the allocation of thermal and electrical load instructions. In general, the primary goal of CHPED is to minimize of the economic costs like fuel costs. With growing concerns about air pollution and other serious environmental issues, the conventional CHPED has already been unable to meet the diversified demands for energy conservation and environmental protection. For this purpose, CHP economic emission dispatch (CHPEED) has been a hot topic since it can take into account environmental protection while pursuing economic benefits [7, 8]. Essentially, a CHPEED problem is to find the optimal heat-power operating point with reasonable fuel costs and emissions, while satisfying a set of various equality and inequality constraints related to heat/electricity demands [9-11]. However, CHPEED poses challenges in terms of computational complexity due to its inherent non-linear, non-convex, and non-smooth characteristic [12], which is hard to solve directly.

## 1.1. Literature review

Research shows that application of multi-objective evolutionary algorithms (MOEAs) is an effective approach for addressing the CHPEED issue and various approaches of this type have been explored in pioneer works [9, 12-17]. In reference [9], the multi-objective particle swarm optimization (MOPSO) is employed to address the CHP stochastic dispatch problem. In reference [12], an enhanced firefly algorithm (EFA) based multi-objective optimization (MOO) method has been put forward to resolve the CHPEED issue. The non-dominated sorting genetic algorithm-II (NSGA-II) is utilized for addressing such problems in reference [13]. In reference [14], a multi-objective line-up competition algorithm (MLCA) has been presented to resolve this issue. In addition, the compromise solution is extracted from the Pareto-optimal solutions by fuzzy decisions. In reference [15], an algorithm based on time-varying acceleration particle swarm optimization (TVA-PSO) is proposed for handling this issue. In reference [16], the non-linear programming (NLP) and the normal boundary intersection (NBI) are used to obtain the Pareto-optimal front (POF) of the problem. Reference [17] summarizes the application of heuristic optimization algorithm in solving CHPEED problems and discusses the implementation of optimization process under different objective functions and constraints. However, to the best of authors' knowledge, no study in the literature has yet identified the best comprise solutions (BCSs) representing decision makers' different preferences with the use of decision analysis in the field of CHPEED, which, to a certain extent, limits the usefulness and practicality of traditional MOEAs-based solution methods.

As a MOO issue, there is no such an optimal solution that enables all objectives to be optimal, and only multiple Pareto-optimal solutions can be obtained in CHPEED [5]. However, it is quite challenging for decision makers to figure out whether a Pareto-optimal solution is a BCS or not from among the noninferior solutions in real-world practice [18]. First, considering the fact that there are a lot of generated Pareto-optimal solutions, the references of decision makers might be different for a specific operation point. Another issue is that for a specific system the preference of the same decision maker may also vary according to changing operation requirements. Therefore, how to identify the BCSs that represent decision makers' different, even conflicting, preferences is a pressing and challenging task for handling CHPEED issues.

### 1.2. Contribution of This Paper

The main contributions of this work are the following three-fold:

(1) A practical CHPEED model: To coordinate the economy and environment protection, a CHPEED model is built with consideration of valve-point loading effects (VPLE) and power transmission losses.

(2) A novel solution approach: A two-stage methodology is proposed for solving the built model. In stage 1, a powerful MOEA, called $\theta$-dominance based evolutionary algorithm ($\theta$-DEA), is employed to seek multiple Pareto-optimal solutions of this model; while stage 2 introduces an integrated decision making (IDM) technique, i.e. fuzzy c-means (FCM)-grey relation projection (GRP), to automatically identify the best comprise solutions reflecting decision makers' different preferences.

(3) The simulation results of three test cases, from simple to complex, prove the effectiveness of our approach, and furthermore, the results demonstrate that our approach is remarkably superior to other state-of-the-art methods.

### 1.3. Organization of This Paper

The rest of this paper is organized as follows. A detailed problem formulation is provided in Section 2, the two-stage solution methodology is shown in Section 3, case studies are investigated in Section 4, and finally conclusions are made in Section 5.

## 2. Problem formulation

## 2.1 Objective functions

### 2.1.1 Fuel costs

The VPLE was originally proposed in reference [19], and then it has been taken into account in recent references [1, 11, 13, 20]. It is the phenomenon that the loss of steam leads to the increase of consumption when a valve is suddenly opened, which makes the unit consumption curve superimposing with pulsating effects [17], as shown in Fig. 1.

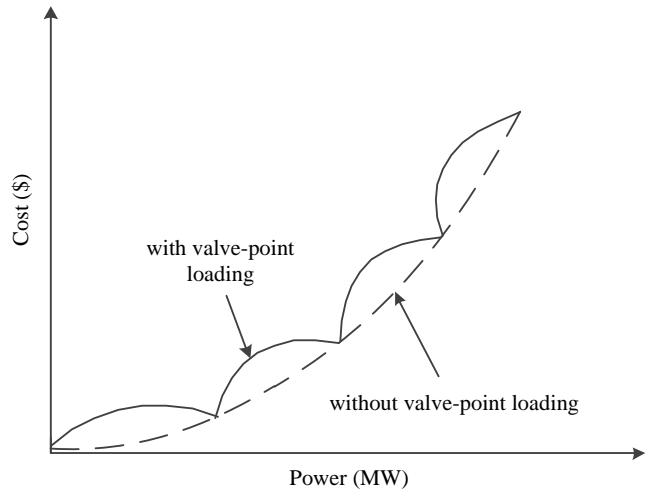

**Fig. 1. Illustration of the VPLE.**

The change function of the consumption characteristics of thermal power units resulting from the VPLE is formulated as $|e_i\sin\{\zeta_i(P_i^{\min} - P_i)\}|$, where $P_i^{\min}$ is the lower limit of the output of unit $i$; $e_i$ and $\zeta_i$ are the cost change coefficients caused by the VPLE of power-only unit $i$, and thereby the fuel cost is modeled as [12]

$$C_{Total} = \sum_{i=1}^{N_p} C_{p,i}(P_i) + \sum_{j=1}^{N_c} C_{c,j}(O_j, H_j) + \sum_{k=1}^{N_h} C_{h,k}(T_k)$$
$$= \sum_{i=1}^{N_p}[a_i + b_i P_i + d_i P_i^2 + |e_i\sin\{\zeta_i(P_i^{\min} - P_i)\}|]$$
$$+ \sum_{j=1}^{N_c}[\alpha_j + \beta_j O_j + \gamma_j O_j^2 + \delta_j H_j + \varepsilon_j H_j^2 + \xi_j O_j H_j] + \sum_{k=1}^{N_h}[\varphi_k + \eta_k T_k + \lambda_k T_k^2]$$
(1)

where $C_{Total}$ is the total fuel cost; $C_{p,i}$, $C_{c,j}$, $C_{h,k}$ are the fuel costs of power-only unit $i$, CHP unit $j$ and heat-only unit $k$, respectively; $a_i$, $b_i$, $d_i$ are the cost coefficients of power-only unit $i$; $\varphi_k$, $\eta_k$, $\lambda_k$ are the cost coefficients of heat-only unit $k$; $\alpha_j$, $\beta_j$, $\gamma_j$, $\delta_j$, $\varepsilon_j$, $\xi_j$ are the cost coefficients of CHP unit $j$; $P_i$, $O_j$ are the generating capacities of power-only unit $i$ and CHP unit $j$, respectively; $H_j$, $T_k$ are the calorific values of CHP unit $j$ and heat-only unit $k$, respectively; $N_p$, $N_h$ and $N_c$ are the numbers of power-only units, heat-only units and CHP units.

### 2.1.2 Gas emission

Since $SO_2$, $NO_x$ and $CO_2$ are particularly harmful to environments, as a major consideration, the gas emission can be modelled as [13]

$$E_{Total} = E_S + E_C \tag{2}$$

$$\begin{aligned} E_S &= \sum_{i=1}^{N_p} E_{p,i}(P_i) + \sum_{j=1}^{N_c} E_{c,j}(O_j) + \sum_{k=1}^{N_h} E_{h,k}(T_k) \\ &= \sum_{i=1}^{N_p} [\mu_i + \kappa_i P_i + \pi_i P_i^2 + \sigma_i e^{(\upsilon_i P_i)}] + \sum_{j=1}^{N_c} \tau_j O_j + \sum_{k=1}^{N_h} \rho_k T_k \end{aligned} \tag{3}$$

$$E_C = \sum_{i=1}^{N_p} \vartheta_i P_i + \sum_{j=1}^{N_c} \psi_j P_j + \sum_{k=1}^{N_h} \varpi_k P_k \tag{4}$$

where $E_{Total}$ is the total emission of polluting gases, $E_S$ is the total emission of $SO_2$ and $NO_x$, $E_C$ is the emission of $CO_2$; $\mu_i$, $\kappa_i$, $\pi_i$, $\sigma_i$, $\upsilon_i$ are the $SO_2$ and $NO_x$ emission coefficients of power-only unit $i$; $\tau_j$ and $\rho_k$ are the $SO_2$ and $NO_x$ emission coefficients of CHP unit $j$ and heat-only unit $k$; $\vartheta_i$, $\psi_j$, $\varpi_k$ are the $CO_2$ emission coefficients of different types of units respectively.

## 2.2 Constraints

### 2.2.1. Power demand constraint

The constraint of power demands is shown as follows [13, 21]:

$$\sum_{i=1}^{N_p} P_i + \sum_{j=1}^{N_c} O_j = P_D + P_L \tag{5}$$

where $P_D$, $P_L$ are the power demand and power transmission loss, in which:

$$P_L = \sum_{i=1}^{N_p}\sum_{j=1}^{N_p} P_i B_{i,j} P_j + \sum_{i=1}^{N_p}\sum_{j=1}^{N_c} P_i B_{i,j} O_j + \sum_{i=1}^{N_c}\sum_{j=1}^{N_c} O_i B_{i,j} O_j + \sum_{i=1}^{N_p} B_{0,i} P_i + \sum_{i=1}^{N_c} B_{0,i} O_i + B_{0,0} \tag{6}$$

In Eq. (6), $B_{i,j}$ is the loss coefficient related to the productions of unit $i$ and $j$; $B_{0,i}$ is the loss coefficient concerned with the production of unit $i$; $B_{0,0}$ is the loss coefficient parameter.

### 2.2.2. Heat demand constraints

The amount of heat required for the system is shown in Eq. (7):

$$\sum_{j=1}^{N_c} H_j + \sum_{k=1}^{N_h} T_k = H_D \tag{7}$$

where $H_D$ denotes the total system heat demand.

### 2.2.3. Capacity limits of each unit

The constraints of each unit are as follows [8, 13]:

$$P_i^{\min} \leq P_i \leq P_i^{\max}, i = 1,...,N_p \tag{8}$$

$$P_j^{\min}(H_j) \leq P_j \leq P_j^{\max}(H_j), j = 1,...,N_c \tag{9}$$

$$H_j^{\min}(P_j) \leq H_j \leq H_j^{\max}(P_j), j = 1,...,N_c \tag{10}$$

$$H_k^{\min} \leq H_k \leq H_k^{\max}, k = 1,...,N_h \tag{11}$$

Eq. (8) denotes the constraint of the power-only units; Eq. (9) and (10) are the constraints of the CHP units; Eq. (11) is the constraint of the heat-only units. For a CHP unit, the generations of heat and power are interrelated and limited each other, and its heat-power feasible operation region (FOR) is illustrated in Fig. 2 [1, 9, 12-15].

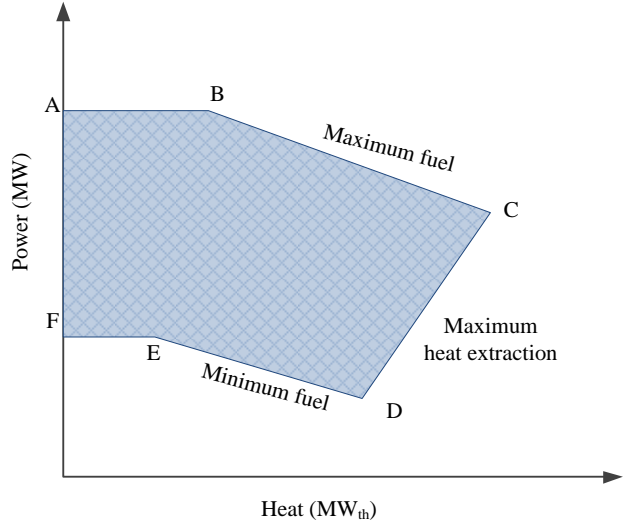

**Fig. 2.** Schematic diagram of FOR.

The closed area surrounded by the curve ABCDEF represents the feasible region of the CHP units. For example, with the increase of heat generation, the unit generating capacity increases in line segment DC, while, in segments BC and DE, the unit generating capacity will decrease.

### 2.2.4. Ramp rate limits

In a dynamic dispatch, each unit must meet the following ramp rate limits [11, 22]:

$$\begin{cases} P_{i,t} - P_{i,(t-1)} \leq \text{UR}_i \\ P_{i,(t-1)} - P_{i,t} \leq \text{DR}_i \end{cases} \quad t = 1, 2, \dots, N_T \quad (12)$$

where $t$ denotes a time interval; $N_T$ is the number of time intervals; $\text{UR}_i$ and $\text{DR}_i$ are respectively the ramp-up and ramp-down rate limits of unit $i$.

## 3. Proposed approach

### 3.1 Solution framework

Different from a mono-objective optimization problem, a MOO produces a set of optimal solutions rather than an optimal solution to coordinate differently weighted or even conflicting objectives. Therefore, it is very suitable to incorporate multiple attribute decision making into the MOO for expressing the preferences of the decision makers such that the BCSs can be identified to better meet the practical needs of system operation [18]. The stage 1 of our approach comprises of the use of $\theta$-DEA [23] to seek the model's Pareto optimal solutions, while the BCSs are determined via the IDM technique FCM-GRP in stage 2. The proposed solution framework is shown in Fig. 3.

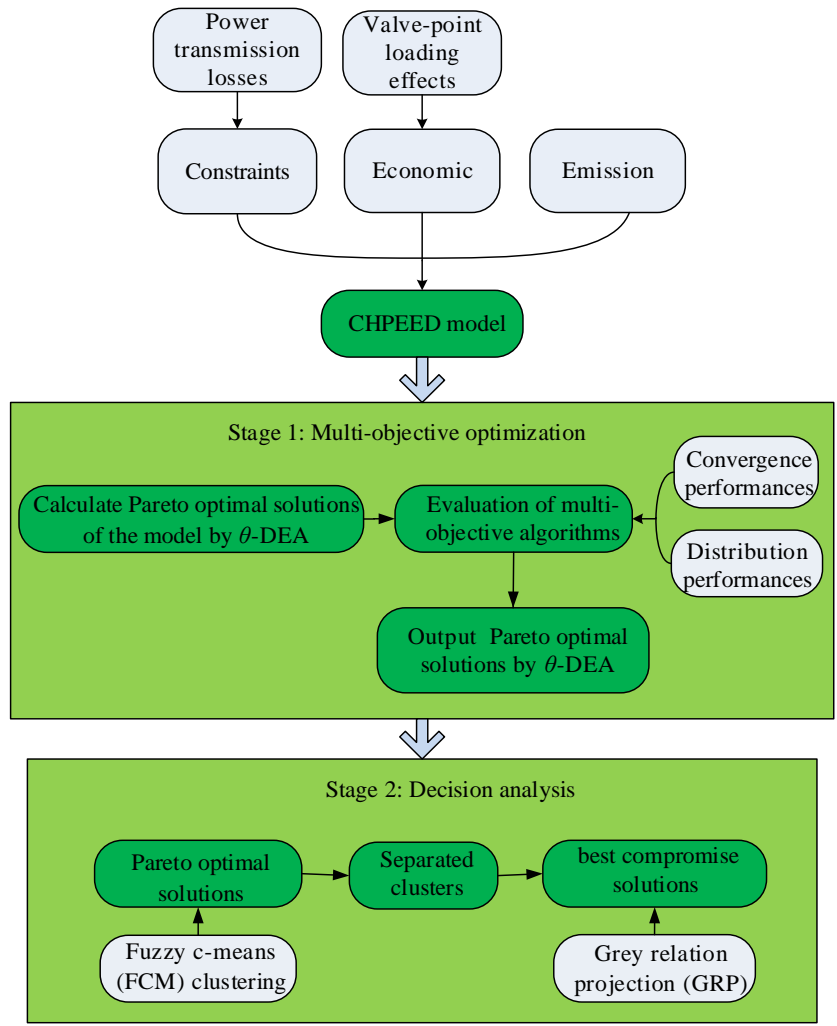

Fig. 3. Proposed two-stage solution framework

### 3.2 Related algorithms

#### 3.2.1 $\theta$-dominance based evolutionary algorithm

Based upon $\theta$ dominance, the $\theta$-DEA proposed by Yuan et al. in 2016 is new powerful algorithm for resolving MOO problems [23].

### 3.2.1.1 θ-Dominance

θ-Dominance is the most important idea in θ-DEA for balancing the trade-off between convergence and diversity. Suppose each solution in a population $S$ is associated with a cluster by a clustering operator, $\tilde{F}(x)=(f_1(x), f_2(x), ..., f_M(x))^T$ is the normalized objective vector of solution $x$, $L$ is a line passing through the origin and the $j$th reference point $\lambda_j$, and **u** is the projection of $\tilde{F}(x)$ on $L$. Let

$$\Upsilon_j(x) = Dis_{j,1}(x) + \theta \times Dis_{j,2}(x), \ j \in \{1, 2, ..., N_{cl}\} \tag{13}$$

where $\theta$ is penalty parameter, $N_{cl}$ is the number of clusters in the population, $Dis_{j,1}(x)$ is the distance between the origin and **u**, $Dis_{j,2}(x)$ is the perpendicular distance between $\tilde{F}(x)$ and $L$. The distances in 2-dimensional objective space is illustrated in Fig. 4.

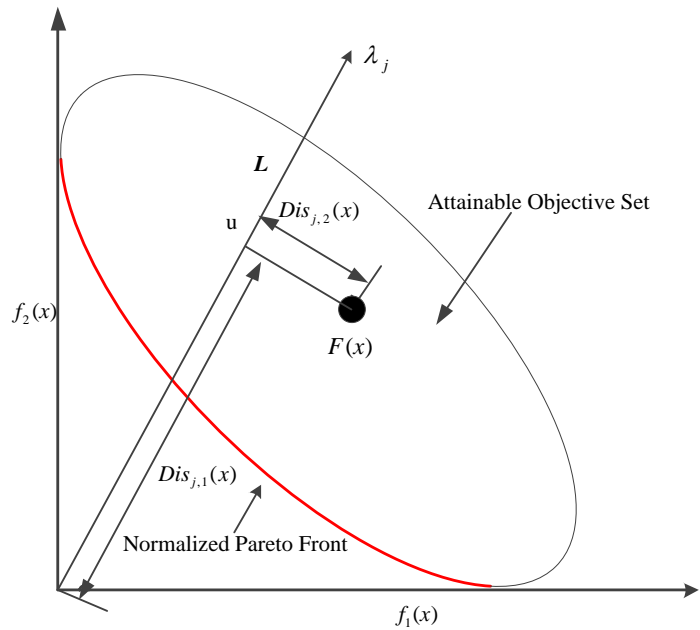

**Fig. 4. Illustration of the distances.**

**Definition 1 [23]:** Given two solutions $x_1$, $x_2$ in the same cluster, $x_1$ is said to θ-dominate $x_2$, denoted by $x_1 \prec_\theta x_2$, and $\Upsilon_j(x_1) < \Upsilon_j(x_2)$, where $j \in \{1, 2, \cdots, N_{cl}\}$

Research demonstrates that through utilizing the $\theta$-Dominance-based fitness evaluation scheme, $\theta$-DEA manages to emphasize convergence and diversity of the algorithm. The more details about the $\theta$-DEA can be found in reference [23].

#### 3.2.1.2 Algorithm flow

The algorithm flow of the $\theta$-DEA is briefly described as follows.

Step 1: Generate $N_d$ reference points $\Lambda = \{\lambda_1, \lambda_2, ..., \lambda_{N_d}\}$.

Step 2: Generate the initial population $A_0$, the ideal initial point $Z^*$ and the worst point $Z^{nad}$.

Step 3: Set the current iteration number $Iter$ to 0.

Step 4: Determine whether the termination criterion is satisfied. If satisfied, output the set of final nondominated solutions $A$ and terminate; otherwise, go to the next step.

Step 5: Generate the offspring population $Q_{Iter}$ via recombined operators, and then obtain a new population $R_{Iter}$ by combining $Q_{Iter}$ with the current population $A_{Iter}$.

Step 6: Obtain the population $S_{iter} = \bigcup F_i$, where $F_i$ is the $i$th Pareto non-dominated levels of $R_{Iter}$.

Step 7: Normalize the population $S_{iter}$, and split it into $N_{cl}$ clusters $CL = \{CL_1, CL_2, ..., CL_{N_{cl}}\}$, where the cluster $CL_j$ is represented by the reference point $\lambda_j$.

Step 8: Classify $S_{Iter}$ into $\theta$-non-domination levels ($F_1^{'}$, $F_2^{'}$, etc.) via the $\theta$-Dominance-based non-dominated sorting.

Step 9: Fill the population slots in $A_{Iter+1}$ with the use of one level at each time.

Step 10: Set $A_{Iter+1}$ to $\phi$, and assign $i$ to 1.

Step 11: Judge whether the condition $|A_{Iter+1}| \cup |F_i^{'}| < N_d$ is met. If met, assign

$A_{Iter+1} \cup F_i^{'}$ to $A_{Iter+1}$ and add 1 to *i*, then repeat this step; otherwise, go to the next step.

Step 12: Randomly sort $F_i^{'}$.

Step 13: assign $A_{Iter+1} \cup F_i^{'}[1:(N_d - |A_{Iter+1}|)]$ to $A_{Iter+1}$.

Step 14: Add 1 to *Iter*, then go to Step 4.

The flowchart of the *θ*-DEA is shown in Fig. 5.

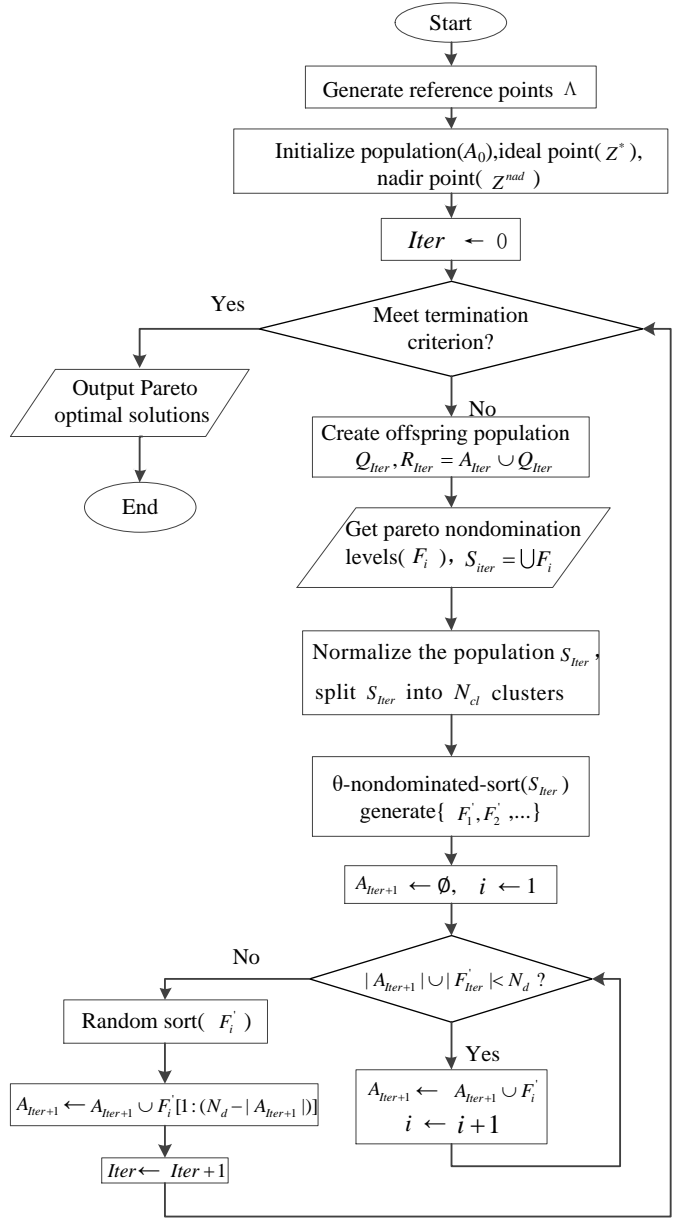

**Fig. 5. Flowchart of the *θ*-DEA.**

### 3.2.2 Fuzzy C-mean Clustering

As a well-known unsupervised clustering algorithm, FCM clustering is based on addressing the following issue [24, 25]:

$$\min \ J(W,U,V) = \sum_{i=1}^{N_d} \sum_{j=1}^{N_{cl}} \mu_{ij}^m \|w_i - v_j\|^2$$
$$s.t. \ \sum_{j=1}^{N_{cl}} \mu_{ij} = 1 \quad (14)$$

where $J$ is a loss function, $W = \{w_1, w_2, \cdots, w_i, \cdots, w_{N_d}\}$ is an input vector with $N_d$ Pareto-optimal solutions to be clustered, $U = \{u_1, u_2, \cdots, u_{N_{cl}}\}$ and $V = \{v_1, v_2, \cdots, v_j, \cdots, v_{N_{cl}}\}$ are respectively the membership degree matrix and cluster centers, $\mu_{ij}$ is a membership degree that represents $x_i$ belongs to cluster $v_j$, $N_{cl}$ is a pre-given number of clusters and $m(m \in [1,\infty])$ is a parameter that determines the degree of fuzziness.

The key idea of FCM is to seek the minimization of the loss function $J$ through repeatedly updating matrix $U$ and cluster $V$ [26]. In order to reflect decision makers' preferences over economy and environmental protection, herein the clustering number $N_{cl}$ is set to 2. By doing so, FCM provides the cluster centroids for the CHPEED issue and separates each Pareto-optimal solution into a proper cluster.

### 3.2.3 Grey Correlation Projection

GRP theory is an effective means to deal with multiple attribute decision making problems with grey information and has been applied in many engineering areas [18, 26]. Since both the objective functions are "benefit-type" evaluation indicators in this problem, the projection $Prj_l$ of a dispatch scheme $l$ onto the ideal scheme is given by

$$Prj_l^{+(-)} = \sum_{i=1}^{N_g} Grc_{l,i}^{+(-)} \times \frac{r_i^2}{\sqrt{\sum_{i=1}^{N_g}(r_i)^2}} \quad (15)$$

where superscript "+"/"−" denote a positive/negative scheme, $N_g$ is the number of indicators used to assess a scheme, $Grc_{l,i}$ is the grey relation coefficient between indicator $i$ and scheme $l$, $r_i$ is the weight of each indicator in a scheme. In this work, for ease of analysis, the weights of the two objectives are set to be equal. The relative projection (*RP*) is defined as follows [26]:

$$RP_l = \frac{Prj_l^+}{Prj_l^+ + Prj_l^-}, \quad 0 \le RP_l \le 1 \quad (16)$$

where $RP_l$ is the relative projection of scheme $l$, which measures the degree of a scheme close to the positive ideal scheme and away from the negative ideal scheme. From Eq. (16), it can be seen that a scheme will be better with a higher relative projection value. Consequently, the solutions with the highest relative projection are considered to be the BCSs.

### 3.4 Evaluation measures of MOEAs

The quantitative evaluation of the performance of MOEAs has recently been attracting concerns. However, there is still no consensus on the evaluation criteria so far in the MOEA community. Generally, a good measure should meet such criteria [18, 23, 27]:

(1) **Criteria 1**——minimum distance: the obtained POF should be as close as possible to the true POF to ensure a good convergence of the algorithm.

(2) **Criteria 2**——uniform distribution and maximum spread: the obtained Pareto-optimal solutions should be uniformly distributed with maximum spread to maintain a good diversity among obtained solutions.

### 3.4.1 Inverted generational distance

As a performance metric, inverted generational distance (IGD) has been widely used for MOO problems in recent studies [23, 27, 28], since it has advantages such as high computational efficiency and wide generality. For a MOEA, let $\Gamma$ and $A$ be the set of targeted points and the set of final nondominated points, this metric as the average Euclidean distance of points in set $\Gamma$ with their nearest members in set $A$ is calculated by [23, 28]

$$IGD(A,\Gamma) = \frac{1}{|\Gamma|} \sum_{i=1}^{|\Gamma|} \min_{n \in A} Edis(\Gamma_i, n) \tag{17}$$

where $Edis(\Gamma_i, n) = \|\Gamma_i - n\|_2$. It should be noted that the set $A$ with smaller values of this measure IGD has better performances about convergence and diversity [28].

### 3.4.2 Spread

The metric Spread [29, 30] is another popular indication for measuring the distribution of a Pareto front obtained in the objective space, as defined by

$$Spread = \frac{Edis_1 + Edis_2 + \sum_{i=1}^{N_d - 1} |Edis_i - \overline{Edis}|}{Edis_1 + Edis_2 + (N_d - 1) \times \overline{Edis}} \tag{18}$$

where $Edis_1$ and $Edis_2$ are the Euclidean distances between extreme solutions and boundary solutions of the set $A$ with $N_d$ members [30]; $\overline{Edis}$ is the average value of all distance $Edis_i$ ($i = 1, \cdots, N_d - 1$). It is not necessary to write that, for two different distributions, this metric with a higher value represents a worse distributions [30].

### 3.5 Solving Process

Fig. 6 illustrates the solving process of our approach, and the details of major procedures are given below.

Step 1: Seek the Pareto-optimal solutions of CHPEED with VPLE and power transmission loss by the $\theta$-DEA.

Step 2: Evaluate the performance of $\theta$-DEA by calculating the quality indicators *IGD* and *Spread*, and then compare with that of other alternatives including MOPSO and NSGA-II.

Step 3: Output Pareto optimal solutions.

Step 4: Initialize the membership degree matrix $U$.

Step 5: Calculate the cluster center $V$.

Step 6: Determine whether the criterion $|J_{cur} - J_{pre}| < \varepsilon$ is satisfied, where $J_{cur}$ and $J_{pre}$ are respectively the current and previous loss function $J$, $\varepsilon$ is a preset threshold. If satisfied, update membership matrix $U$ and go to Step 5; otherwise, output separated clusters.

Step 7: Base on the obtained clusters, build an initial decision matrix.

Step 8: Standardize the decision matrix in Step 7.

Step 9: Calculate grey relation coefficient $Grc$ between indicator and scheme.

Step 10: Calculate priority $Prj_l$ of scheme $l$ onto the ideal schemes according to Eq. (15).

Step 11: Calculate the relative projection $RP_l$ of scheme $l$ according to Eq. (16).

Step 12: Output the BCSs with the highest $RP$ values.

Step 13: Compare with the existing results in the previous literature.

Step 14: Output the optimal dispatch schemes.

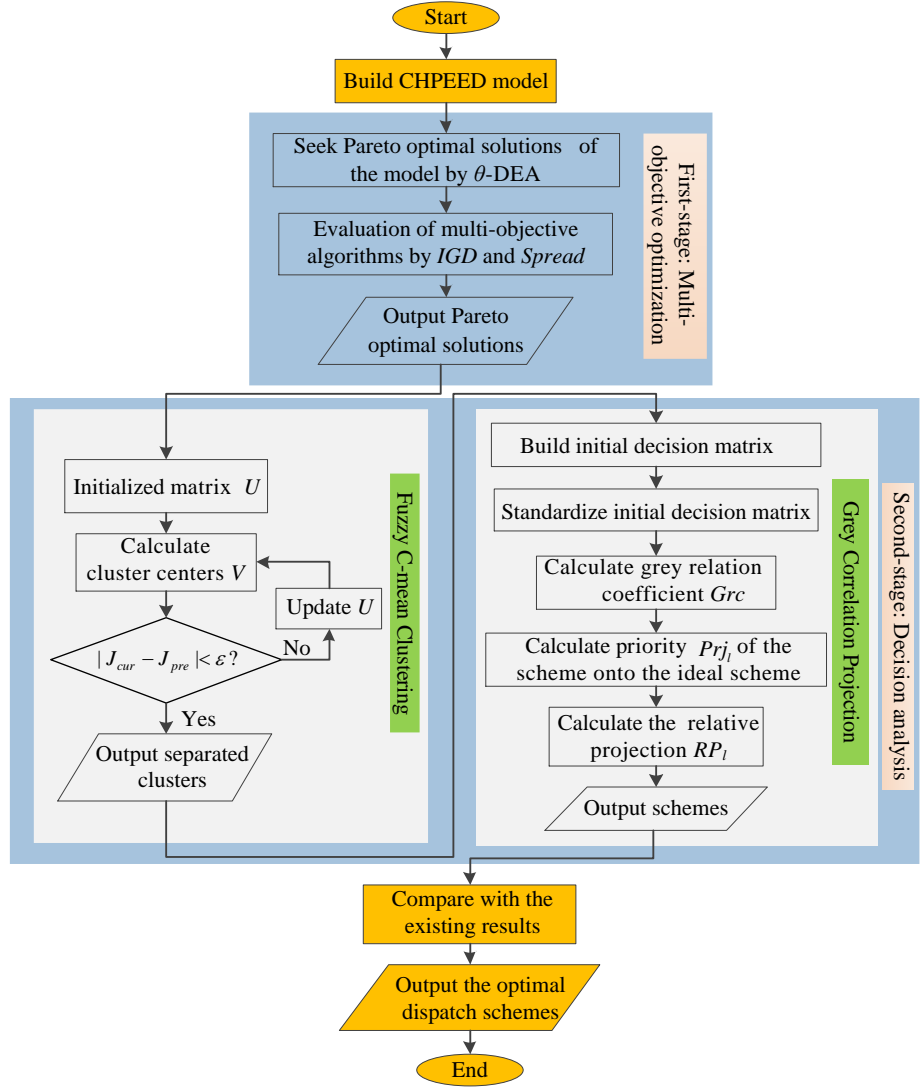

**Fig. 6. Solving process of the proposed approach**

## 4. Case Studies

In order to test the performance of our method, three test cases, from simple to complex, have been set up in this section. Specially, Case 1 only takes into account the VPLE, while Cases 2 and 3 simultaneously consider both the VPLE and power transmission loss, wherein Case 3 is a dynamic dispatch problem. To further demonstrate the strong competitiveness, an extensive comparison of our approach with other advanced algorithms has been made. For ease of comparison, both the maximum iteration number and population size are set to 100 in all the algorithms. All simulations are implemented on a PC platform with 2 Intel Core dual core CPUs (2.4 GHz) and 6 GB RAM.

**4.1. Case 1**

This test case, originally proposed by Gou et al. [31], includes one power-only unit (unit 1), three CHP units (units 2-4) and one heat-only unit (unit 5). Based on references [13, 31, 32], the specific model parameters used in this test case are described in detail in the Appendix A.1.

For the purpose of properly evaluating the performance of the $\theta$-DEA, two commonly used MOEAs, i.e. MOPSO [9] and NSGA-II [13], are employed to derive Pareto-optimal solutions.

The POFs using these algorithms in this Case are shown in Fig. 7.

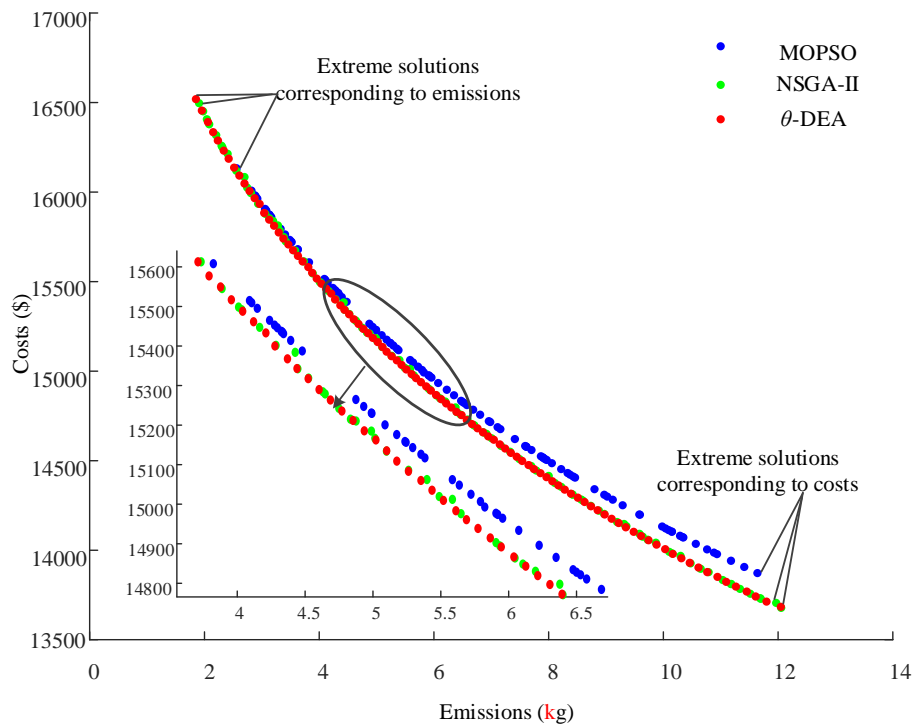

**Fig. 7. Distribution of POFs in Case 1.**

Form Fig. 7, it is clear that the $\theta$-DEA is capable of generating nearly complete and uniform Pareto-optimal solutions with well-distribution. Thereafter, a conclusion can be drawn that the economy and environment protection in CHPEED issues can be effectively coordinated via the proposed approach. Furthermore, Fig. 7 also indicates that the optimization ability of the $\theta$-DEA is superior to that of the NSGA-II and MOPSO algorithms, embodying that the POF obtained by the $\theta$-DEA dominates the fronts of the two comparison algorithms in most cases.

The POF obtained by the $\theta$-DEA before and after clustering are respectively shown in Figs. 8 and 9.

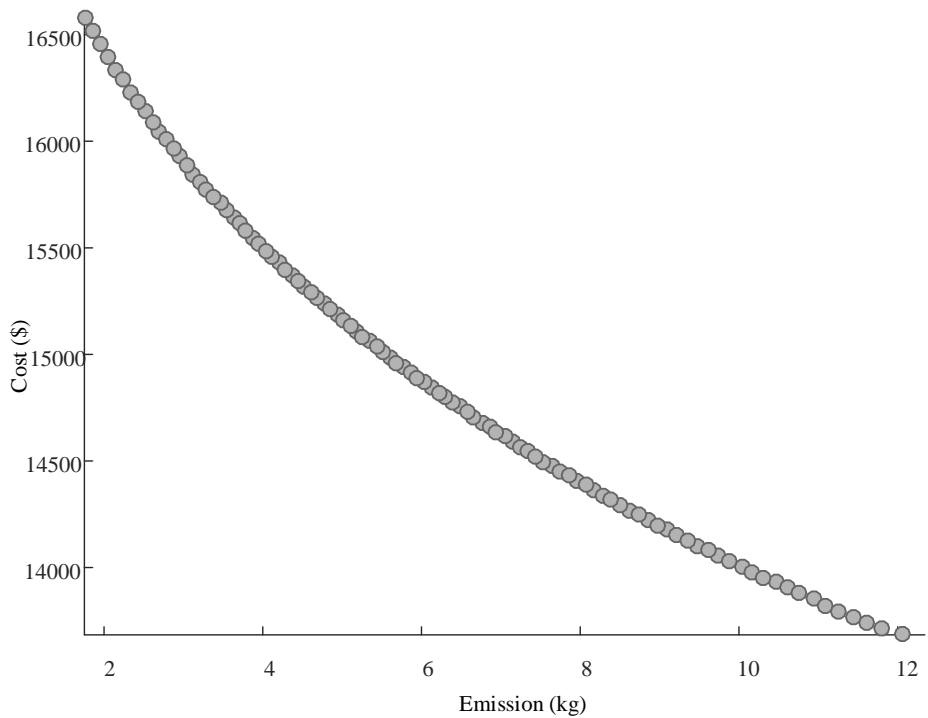

**Fig. 8. POF using the $\theta$-DEA in Case 1 before clustering.**

As is shown in Fig. 8, the fuel costs decrease monotonically with the increase of gas emission, and vice versa. The reason for this phenomenon is that the two objective functions conflict to each other, and there is a trade-off between fuel cost and gas emission. This trade-off makes the determination of BCSs reflecting decision makers' different preferences from all Pareto-optimal solutions even more complex, and thus new decision analysis approaches are required to balance multiple conflicting objectives. For this purpose, the new IDM approach, called FCM-GRP, is employed in this work.

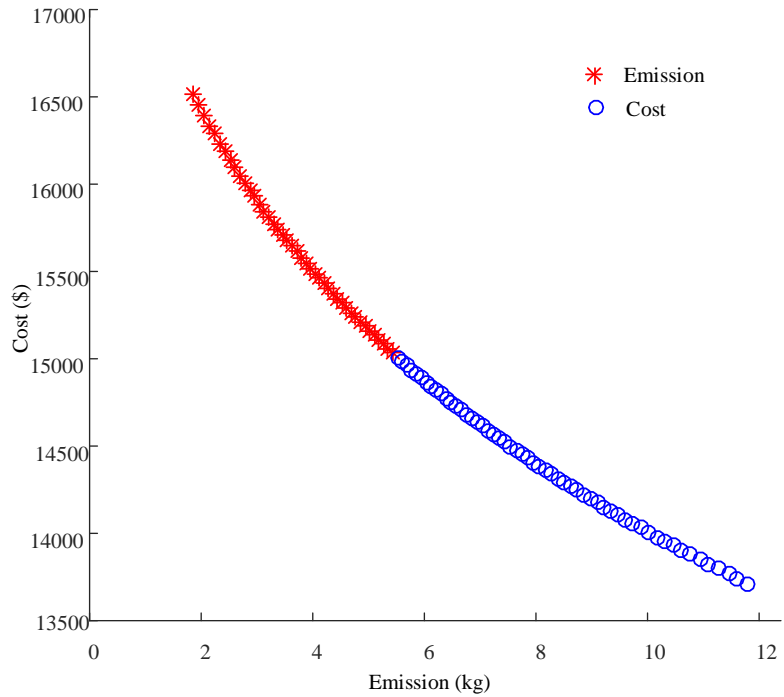

**Fig. 9.** POF using the $\theta$-DEA in Case 1 after clustering.

As shown in Fig. 9, the Pareto-optimal solutions are separated into two groups through the use of FCM clustering. And then, two BCSs extracted from the separated clusters using FCM-GRP in this work and the BCSs obtained by NSGA-II and strength Pareto evolutionary algorithm 2 (SPEA 2) in reference [13] are shown in Table 1.

**Table 1**

Results obtained from our approach and other methods in Case 1 (MW).

| Items | Proposed apprach | | NSGA-II [13] | SPEA 2 [13] |
|---|---|---|---|---|
| | BCS 1 | BCS 2 | | |
| Cost ($) | 14504.2 | 15137.3 | 15008.7 | 14964.3 |
| Emission (kg) | 7.5 | 5.1 | 6.1 | 6.4 |
| $P_1$ (MW) | 105.6 | 85.6 | 93.9 | 96.5 |
| $P_2$ (MW) | 61.7 | 80.1 | 72.8 | 71.2 |
| $P_3$ (MW) | 27.8 | 29.3 | 43.3 | 44.5 |
| $P_4$ (MW) | 104.9 | 105.0 | 89.9 | 87.8 |
| $H_2$ (MWth) | 76.4 | 76.5 | 84.9 | 84.8 |
| $H_3$ (MWth) | 39.5 | 39.3 | 22.6 | 10.2 |
| $H_4$ (MWth) | $5.3*10^{-3}$ | $0.4*10^{-3}$ | 2.6 | 17.9 |
| $H_5$ (MWth) | 34.1 | 34.2 | 39.8 | 37.1 |
| $P_D$ | 300 | 300 | 300 | 300 |
| $H_D$ | 150 | 150 | 150 | 150 |

From Table 1 it can be observed that the optimal results of different algorithms are different from each other. Compared with the BCSs of NSGA-II and SPEA 2, the cost of BCS 1 is respectively decreased by $ 504.5 and $ 460.1, while its emission is increased by 1.4 kg and 1.1 kg; regarding BCS 2, the cost is increased by $ 128.6 and $ 173.0, while the emission is reduced by 1.0 kg and 1.3 kg, respectively. Therefore, the conclusion can be drawn that FCM-GRP is an effective tool to automatically determine the BCSs for the CHPEED problem, which helps to provide more realistic options representing decision makers' different references.

**4.2 Case 2**

In order to further examine the effectiveness of our approach, it is performed in a more complex case. This test case, originally proposed in [31], is widely used in previous works [13, 33-35]. It consists of four power-only units (units 1-4), two CHP units (units 5 and 6) and a heat-only unit (unit 7). The detailed model parameters of this test system are given in the Appendix A.2.

To improve the practicality of the model, different from Case 1, the VPLE and power transmission losses are simultaneously considered in this case. Considering the randomness of MOEAs to optimal results [36], all the employed algorithms are independently performed 30 times, and the comparison results of the evaluation metrics using these algorithms are shown in Table 2.

**Table 2**

Comparison of the evaluation metrics using the algorithms in Case 2.

| Algorithms | Metrics | Average value | Best value | Worst value |
|---|---|---|---|---|
| $\theta$-DEA | *IGD* | 10293.84 | 10801.25 | 10134.20 |
|  | *Spread* | 0.89 | 0.93 | 0.85 |
| MOPSO | *IGD* | 10779.00 | 11375.52 | 10457.21 |
|  | *Spread* | 0.97 | 0.99 | 0.95 |
| NSGA-II | *IGD* | 10321.51 | 10989.91 | 10122.01 |
|  | *Spread* | 0.92 | 0.94 | 0.90 |

From Table 2, it can be found that the metrics IGD and Spread of the $\theta$-DEA are superior to those of the MOPSO and NSGAII. This results indicates that the $\theta$-DEA has better convergence and distribution performances than the alternatives in this case.

The POFs obtained by the three algorithms in Case 2 are illustrated in Fig. 10.

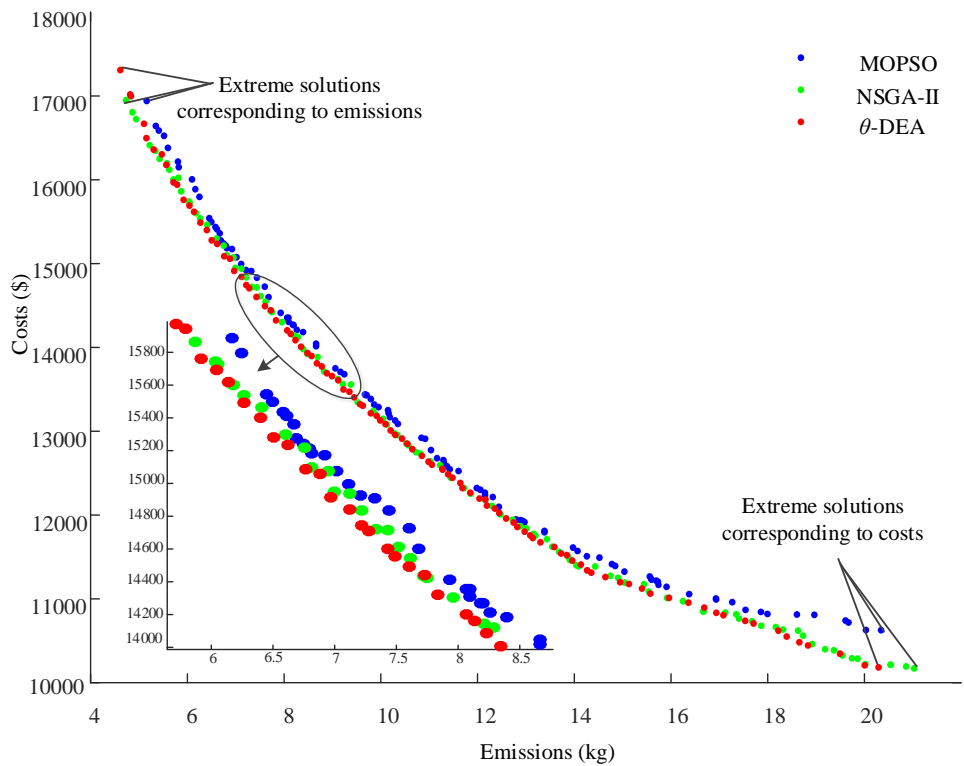

**Fig. 10. Distribution of the POFs in Case 2.**

As shown in Fig. 10, the $\theta$-DEA yields a POF with well-distributed and well-separated solutions. Moreover, the POF of the $\theta$-DEA dominates that of the other algorithms in general. This evidence further validates the effectiveness and superiority of the $\theta$-DEA in finding multiple Pareto-optimal solutions.

Figs. 11 and 12 shows the extreme solutions and the corresponding power transmission loss obtained by $\theta$-DEA, together with the existing results in literature using various algorithms including real coded genetic algorithm (RCGA) [37], particle swarm optimization (PSO) [38], evolutionary programming (EP) [39], artificial-immune system optimization (AIS) [40], differential evolution (DE) [41], bee colony optimization (BCO) [35] and MOPSO.

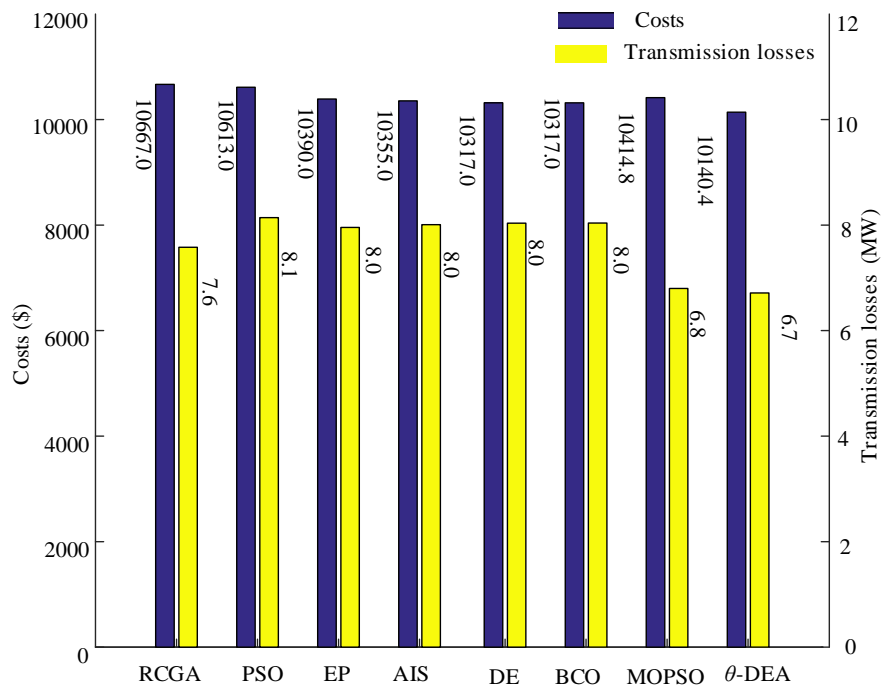

**Fig. 11. Cost extreme solution and transmission losses in Case 2.**

As can be seen in Fig. 11, the cost extreme solution and power transmission loss of the $\theta$-DEA are obviously outperform those of other optimization algorithms.

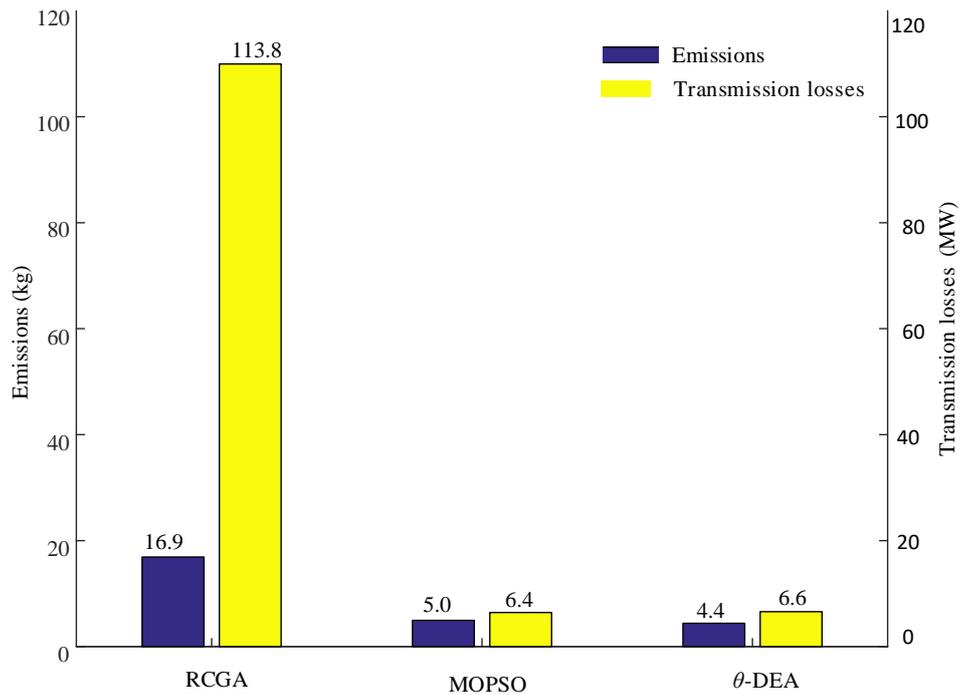

**Fig. 12. Emission extreme solution and transmission losses in Case 2.**

From Fig. 12, it can be observed that the emission extreme solution of the $\theta$-DEA is superior to that of the others without loss of transmission losses. Based on the facts shown in Figs. 11 and 12, we can see that the $\theta$-DEA still has powerful global optimization ability even though being used for single-objective optimization. More importantly, the results show that the extreme solutions of the $\theta$-DEA are better than those of other MOEAs such as RCGA and MOPSO.

The unit outputs of the extreme solutions are demonstrated in Figs. 13 and 14.

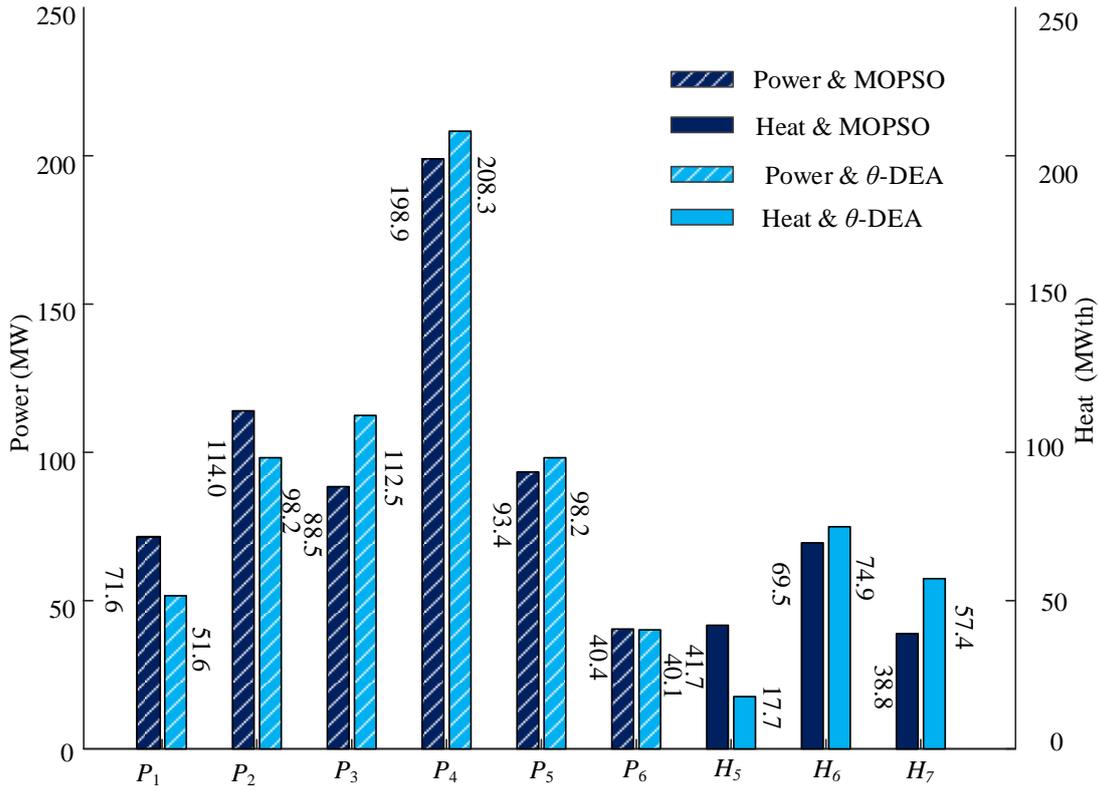

**Fig. 13. Unit outputs of the cost extreme solution in Case 2.**

Fig. 13 indicates that in terms of the cost extreme solution, the power-only units operate near the maximum power points to balance supply and demand satisfying their operating constraints, while the CHP units only play a supporting role. The reason causing such phenomenon is that the power-only units (especially unit 4) yields far less operating costs than other types of units through the analysis of the model parameters in this case, as shown in the Appendix A.2.

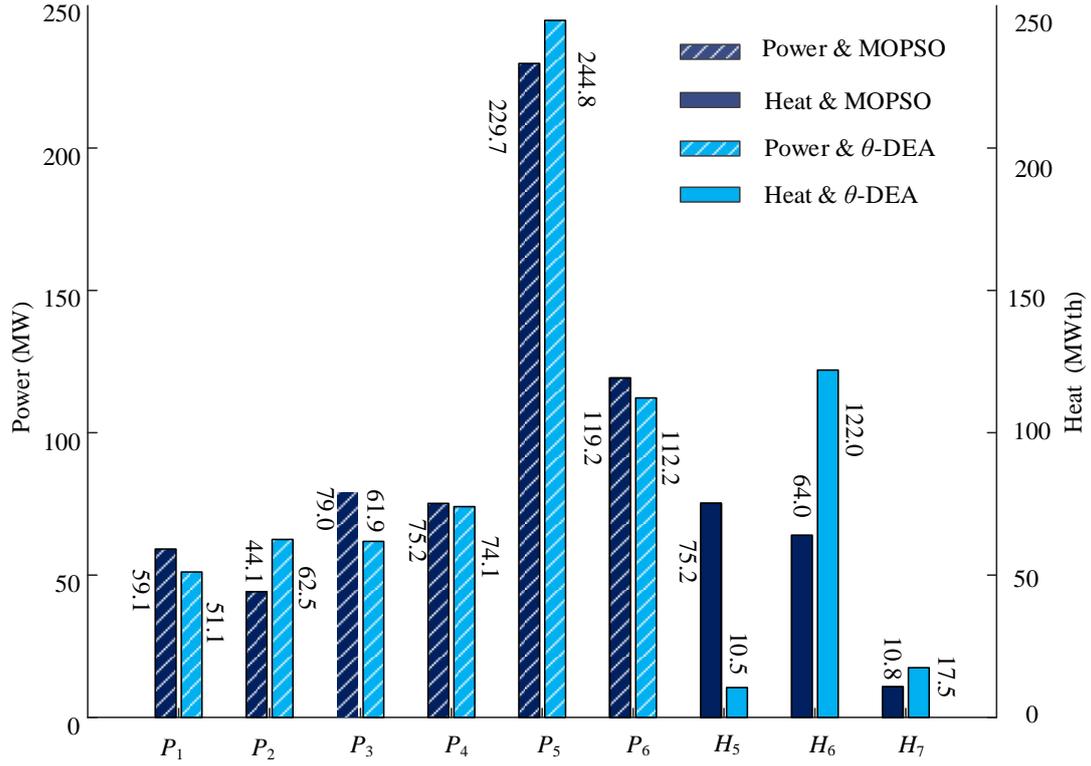

**Fig. 14. Unit outputs of the emission extreme solution in Case 2.**

As seen from Fig. 14, as far as the emission extreme solution is concerned, the CHP units generate electricity as much as possible within their feasible operation regions, while other units generate relatively less electricity. This is due to the fact that emissions of the CHP units are much less than those of others for generating the same amount of electricity which can be seen from the model parameters used in this test case.

To reasonably assess the effectiveness of our approach, a comparative analysis between the obtained BCSs and the existing results of other advanced algorithms such as NSGA-II, MLCA and NBI in literature. The results are shown in Table 3.

**Table 3**

Results obtained from the proposed approach and other methods in Case 2 (MW).

| Items | NSGA-II [13] | MLCA [14] | NBI [33] | Proposed apprach | |
| --- | --- | --- | --- | --- | --- |
| | | | | BCS 1 | BCS 2 |
| Cost ($) | 13433.2 | 12451.4 | 12424.1 | 12196.3 | 13282.9 |
| Emission (t) | 25.8 | 11.1 | 10.9 | 12.0 | 9.7 |
| $P_L$ (MW) | 110.0 | \ | 7.6 | 6.1 | 6.5 |
| $P_1$ (MW) | 73.6 | 62.9 | 68.9 | 69.7 | 64.1 |
| $P_2$ (MW) | 106.9 | 98.5 | 84.1 | 86.6 | 77.7 |
| $P_3$ (MW) | 119.0 | 100.0 | 96.0 | 101.0 | 91.1 |
| $P_4$ (MW) | 163.6 | 105.9 | 117.1 | 123.7 | 113.8 |
| $P_5$ (MW) | 188.4 | 193.0 | 201.4 | 183.1 | 218.0 |
| $P_6$ (MW) | 58.5 | 40.4 | 40.0 | 42.0 | 41.8 |
| $H_5$ (MWth) | 26.8 | 4.7 | 2.5 | 0.7 | 5.5 |
| $H_6$ (MWth) | 74.0 | 75.3 | 75.0 | 76.3 | 74.7 |
| $H_7$ (MWth) | 49.2 | 70.0 | 72.5 | 73.0 | 69.8 |

Table 3 demonstrates that our method has the ability to provide multiple BCSs for decision makers to choose according their own preferences. Specially speaking, compared with the NSGA-II, MLCA and NBI algorithms, the cost is reduced by $ 1236.9, $ 255.1 and $ 227.8 in BCS 1, while the emission is respectively reduced by 16.1 kg, 1.4 kg, and 1.2 kg in BCS 2. At the same time, the power transmission losses of the two BCSs are superior to those of others. As a result, the effectiveness and superiority of the presented method for solving CHPEED issues have been further verified based on these results.

**4.3 Case 3**

This testing system, first studied in [42], is made up of ten units with non-smooth fuel cost and emission level functions, and the model parameters of this system are given in reference [16]. Unlike the static dispatch in the former two cases, a dynamic dispatch is performed to obtain proper power allocations which incorporates the unit generating capability information given by unit ramping constraints, together with the economic and environmental considerations, in this case. The time horizon is herein divided into 24 intervals to satisfy the demand of the system. The power output of each unit in BCSs 1 and 2 during each time interval are shown in Figs. 15 and 16.

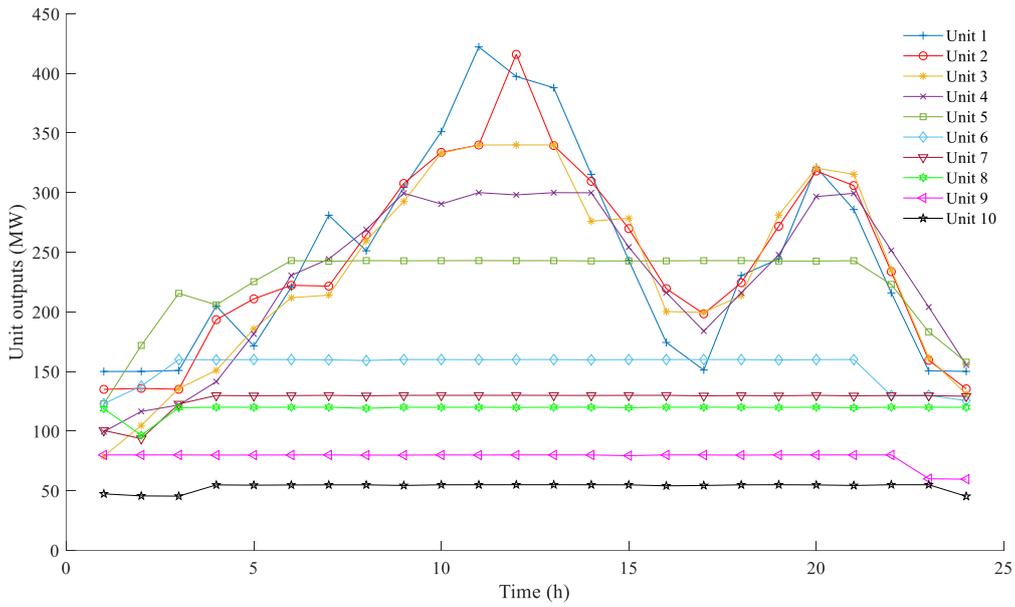

**Fig. 15. Unit outputs in BCS 1 of the proposed approach in Case 3.**

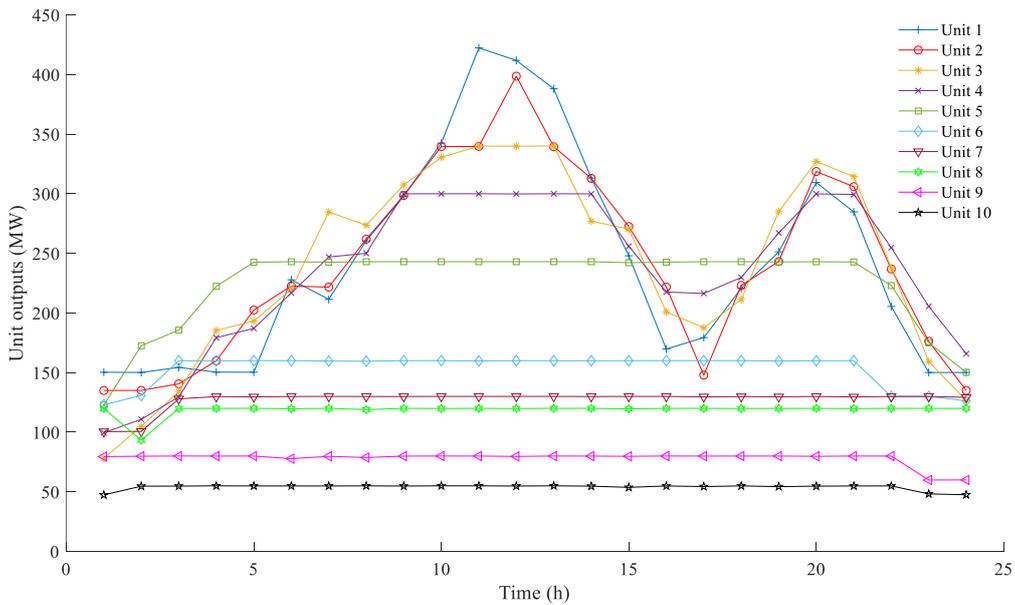

**Fig. 16. Unit outputs in BCS 2 of the proposed approach in Case 3.**

As shown in Figs. 15 and 16, the unit outputs in the BCSs not only real-time balancing the load demands while satisfying all operational constraints including the ramp rate limits. From the results, a conclusion can be reached that our approach is also applicable for solving the dynamic dispatch problem.

To evaluate the performance of our approach for addressing this dynamic dispatch issue, the comparisons of results from the proposed algorithm and other methods are further carried out, in which the results of NSGA-II-based method and scalar optimization (SO)-based algorithm are directly obtained from references [42] and [16].

Figs. 17 and 18 show the costs and emissions of the BCSs in each time period using the proposed approach and two alternatives.

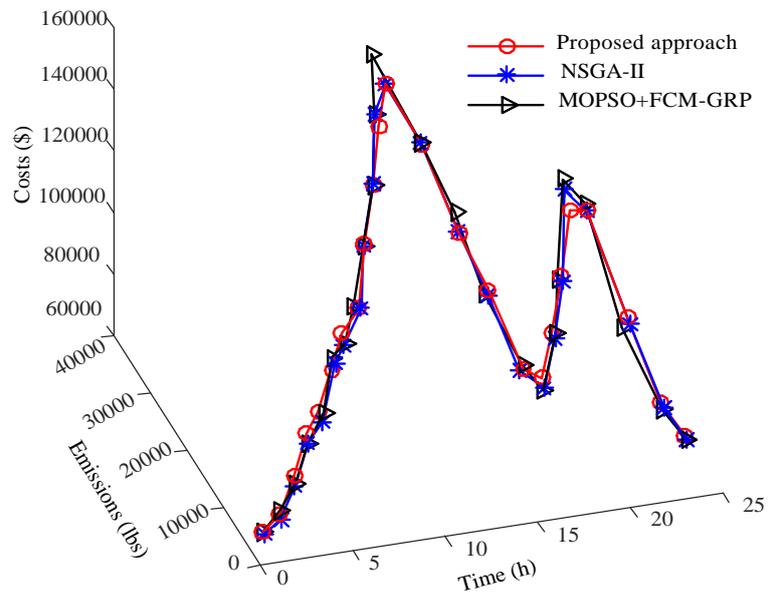

Fig. 17. Costs and emissions of BCS 1 at different intervals.

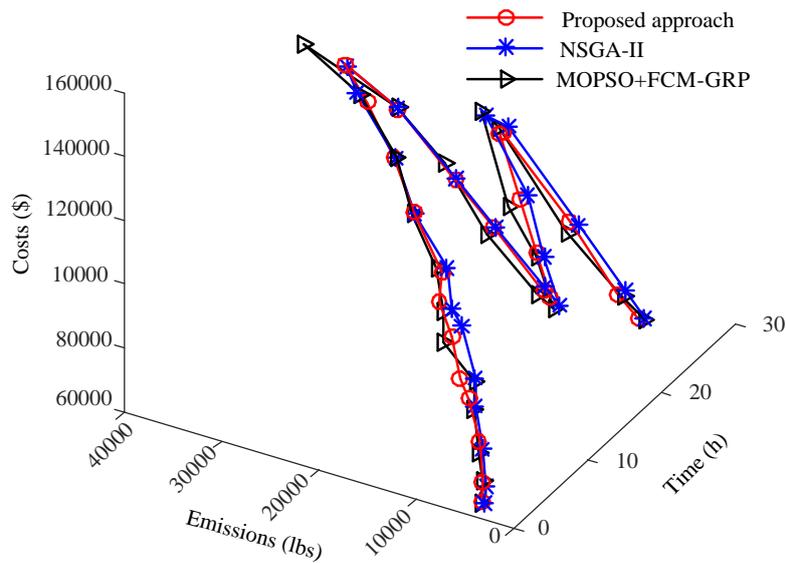

Fig. 18. Costs and emissions of BCS 2 at different intervals.

The total costs and emissions of the BCSs obtained by different algorithms are shown in Fig. 19.

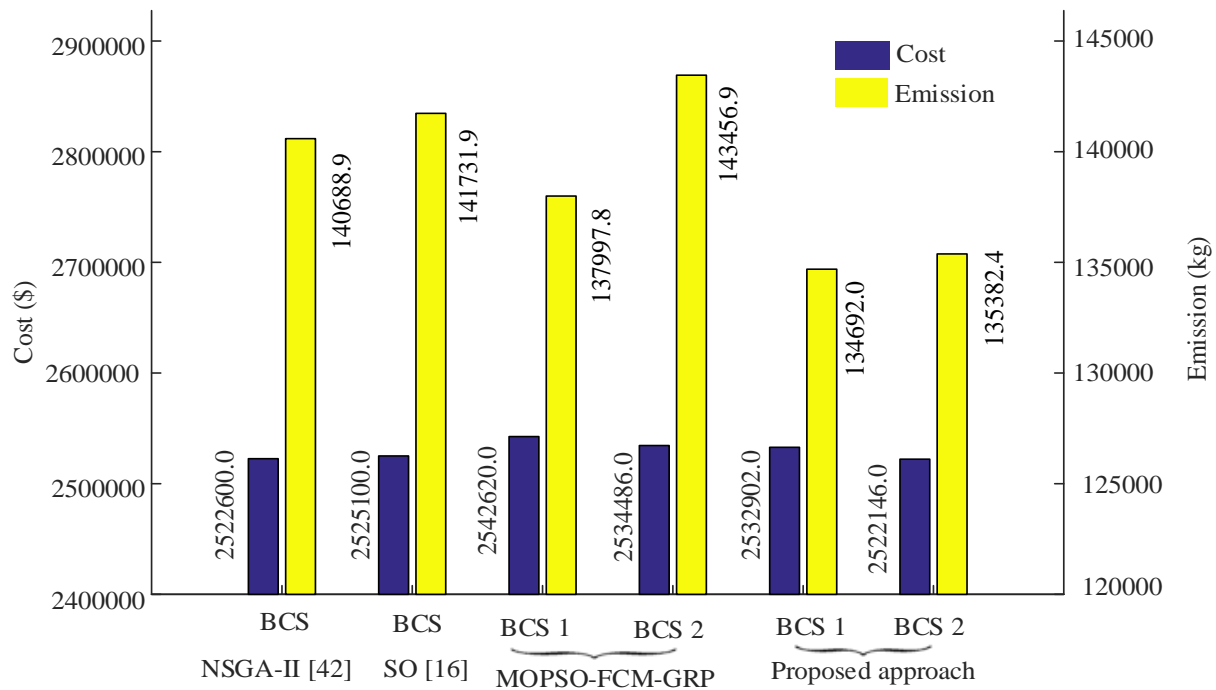

**Fig. 19. Costs and emissions of the BCSs obtained by different algorithms.**

As can be seen from Fig. 19, compared with conventional methods such as NSGA-II [42] and SO [16], the proposed two-stage solution methodology incorporating decision analysis yields multiple BCSs to better reflecting decision maker' different preferences, which will help to improve the practicality of the dispatching strategies obtained through optimization. Furthermore, with regard to specific results, the emission in BCS 1 of our propose approach is less than that of other methods, and the cost in BCS 1 is slightly more than those of NSGA-II (0.41%) and SO (0.31%) but less than the integration method combining MOPSO with FCM-GRP; while the cost and emission in BCS 2 of the propose approach are superior to those of others. Therefore, the conclusion that can be drawn on the basis of the above facts is that the proposed approach manages to solve the dynamic economic emission dispatch problem. Moreover, the superiority of our approach relative to other famous algorithms has also been verified.

## 5. Conclusion

To better meet the diversified demands for energy conservation and environmental protection, a practical CHPEED model is built by taking into account VPLE and power transmission losses. To solve this model, a two-stage solution methodology is proposed by combining multi-objective optimization using the $\theta$-DEA with an integrated decision making technique FCM-GRP. The simulation results of the test cases considered in this paper, from simple to complex, reveal that the proposed approach manages to yield multiple complete and well-distributed Pareto-optimal solutions, but also can automatically identify the BCSs that represent decision makers' different preferences using via the FCM-GRP decision analysis. By doing so, the economy and environment protection in CHPEED issues can be effectively coordinated. Taking Case 1 as an example, the cost and emission of BCS 1 are respectively $ 14504.2 and 7.5 kg, while those of BCS 2 are $ 15137.3 and 5.1 kg. Compared with the results of the NSGA-II and SPEA 2 in [13], the cost of BCS 1 is respectively decreased by $ 504.5 and $ 460.1, while the emission of BCS 2 is reduced by 1.0 kg and 1.3 kg. Therefore, the proposed methodology provides an innovative tool for addressing the CHPEED problem, which manages to provide more realistic options representing decision makers' different references. In addition, the results of this work suggest significant practical implications for determining the best compromise solutions from all Pareto-optimal solutions, which is especially helpful to meet the diverse needs under changing operating conditions of a CHP system.

Our future work will focus on extending this study to extensive potential applications in the optimal operation and control for a smart integrated energy system. In addition, more realistic modeling techniques such as load and renewable generation uncertainties [43], and energy storage units [44] will be incorporated to improve the practicality of the proposed method.

**Acknowledgements**

This work is supported by the U.S. Department of Energy (DOE)'s Office of Electricity Delivery and Energy Reliability – Advanced Grid Modeling (AGM) Program and the China Scholarship Council (CSC) under Grant No. 201608220144.

**Appendix A.**

**A.1. Test system 1**

**A.1.1 Cost and emission function of each unit**

$P_D = 300$ MW, $H_D = 150$ MWth

**(a) Power-only units**

$C_1(P_1) = 254.8863 + 7.6997 P_1 + 0.00172 P_1^2 + 0.000115 P_1^3$ $

$35 \leq P_1 \leq 135$ MW

$E_1(P_1) = 10^{-4} \times (4.091 - 5.554 P_1 + 6.490 P_1^2) + 2 \times 10^{-4} \times \exp(0.02857 P_1)$ kg

**(b) CHP units**

$C_2(P_2, H_2) = 1250 + 36 P_2 + 0.0435 P_2^2 + 0.6 H_2 + 0.027 H_2^2 + 0.011 P_2 H_2$ $

$E_2(P_2, H_2) = 0.00165 P_2$ kg

$C_3(P_3, H_3) = 2650 + 34.5 P_3 + 0.1035 P_3^2 + 2.203 H_3 + 0.025 H_3^2 + 0.051 P_3 H_3$ $

$E_5(P_5, H_5) = 0.0022 P_3$ kg

$C_4(P_4, H_4) = 1565 + 20P_4 + 0.072P_4^2 + 2.3H_4 + 0.02H_4^2 + 0.04P_4H_4$ \$

$E_4(P_4, H_4) = 0.0011 P_4$ kg

**(c) Heat-only units**

$C_5(H_5) = 950 + 2.0109H_5 + 0.038H_5^2$ \$

$0 \leq H_5 \leq 60$ MWth

$E_5(H_5) = 0.0017 H_5$ kg

**A.1.2 Heat-power feasible operation region (FOR) of cogeneration units**

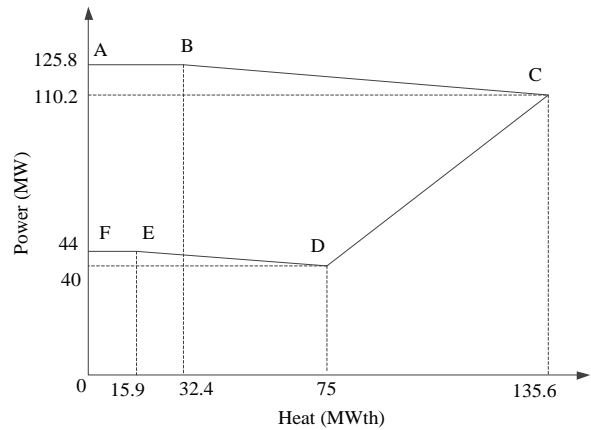

**Fig. 20. FOR of cogeneration unit 1 in test system 1.**

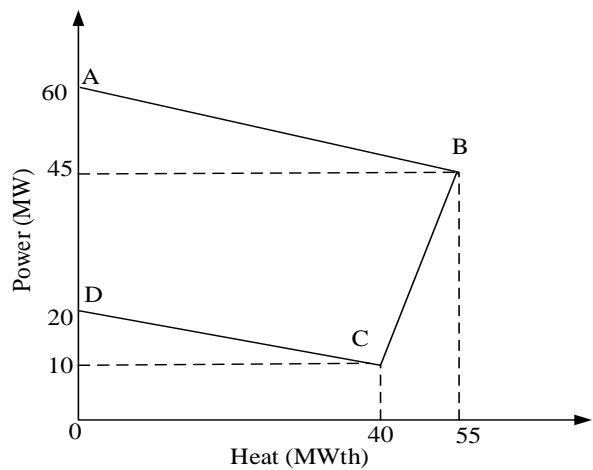

**Fig. 21. FOR of cogeneration unit 2 in test system 1.**

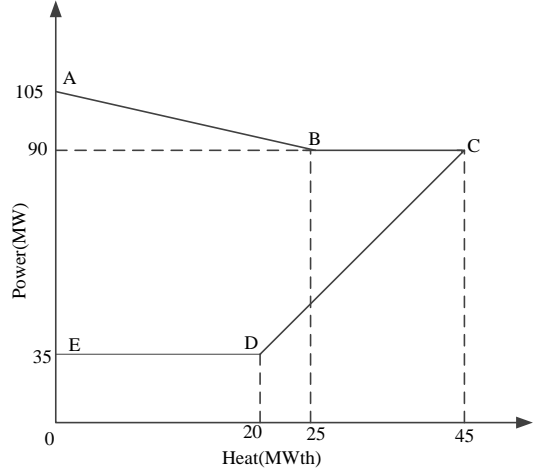

**Fig. 22. FOR of cogeneration unit 3 in test system 1.**

## A.2. Test system 2

### A.2.1 Cost and emission function of each unit

$P_D=600$ MW, $H_D=150$ MWth

**(a) Power-only units**

$C_1(P_1)=25+2P_1+0.008P_1^2+|100\sin\{0.042(P_1^{min}-P_1)\}|$ \$

$10 \leq P_1 \leq 75$ MW

$E_1(P_1)=10^{-4}\times(4.091-5.554P_1+6.490P_1^2)+2\times10^{-4}\times\exp(0.02857P_1)$ kg

$C_2(P_2)=60+1.8P_2+0.003P_2^2+|140\sin\{0.04(P_2^{min}-P_2)\}|$ \$

$20 \leq P_2 \leq 125$ MW

$E_2(P_2)=10^{-4}\times(2.543-6.047P_2+5.638P_2^2)+5\times10^{-4}\times\exp(0.03333P_2)$ kg

$C_3(P_3)=100+2.1P_3+0.0012P_3^2+|160\sin\{0.038(P_3^{min}-P_3)\}|$ \$

$30 \leq P_3 \leq 175$ MW

$E_3(P_3)=10^{-4}\times(4.285-5.094P_3+4.586P_3^2)+1\times10^{-6}\times\exp(0.08P_3)$ kg

$C_4(P_4)=120+2P_4+0.001P_4^2+|180\sin\{0.037(P_4^{min}-P_4)\}|$ \$

$40 \leq P_4 \leq 250$ MW

$$E_4(P_4) = 10^{-4} \times (5.326 - 3.550P_4 + 3.370P_4^2) + 2 \times 10^{-3} \times \exp(0.02P_4) \text{ kg}$$

**(b) CHP units**

$$C_5(P_5, H_5) = 2650 + 14.5P_5 + 0.0345P_5^2 + 4.2H_5 + 0.03H_5^2 + 0.031P_5H_5 \text{ \$}$$

$$E_5(P_5, H_5) = 0.00165P_5 \text{ kg}$$

$$C_6(P_6, H_6) = 1250 + 36P_6 + 0.0435P_6^2 + 0.6H_6 + 0.027H_6^2 + 0.011P_6H_6 \text{ \$}$$

$$E_6(P_6, H_6) = 0.00165P_6 \text{ kg}$$

**(c) Heat-only units**

$$C_7(H_7) = 950 + 2.0109H_7 + 0.038H_7^2 \text{ \$}$$

$$0 \leq H_7 \leq 2695.2 \text{ MWth}$$

$$E_7(H_7) = 0.0018H_7 \text{ kg}$$

**(d) Network loss coefficients**

$$B = \begin{bmatrix} 49 & 14 & 15 & 15 & 20 & 25 \\ 14 & 45 & 16 & 20 & 18 & 19 \\ 15 & 16 & 39 & 10 & 12 & 15 \\ 15 & 20 & 10 & 40 & 14 & 11 \\ 20 & 18 & 12 & 14 & 35 & 17 \\ 25 & 19 & 15 & 11 & 17 & 39 \end{bmatrix} \times 10^{-6}$$

$$B_0 = [-0.3908 \quad -0.1297 \quad 0.7047 \quad 0.0591 \quad 0.2161 \quad -0.6635] \times 10^{-3}$$

$$B_{00} = 0.056$$

**A.2.2 Heat-power feasible operation region of cogeneration units**

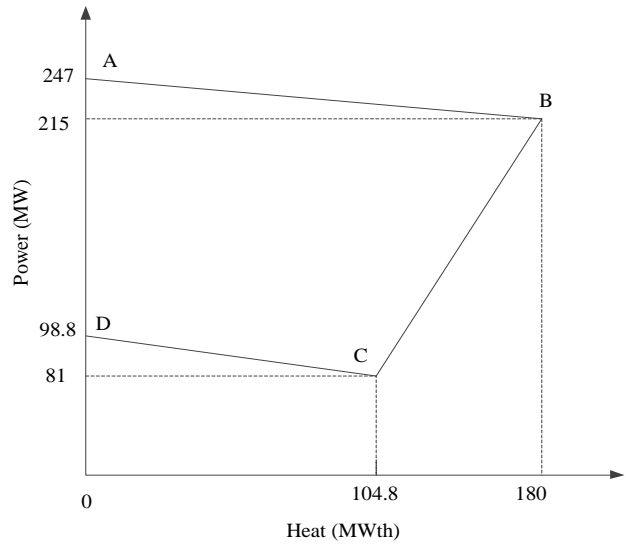

**Fig. 23. FOR of cogeneration unit 1 in test system 2.**

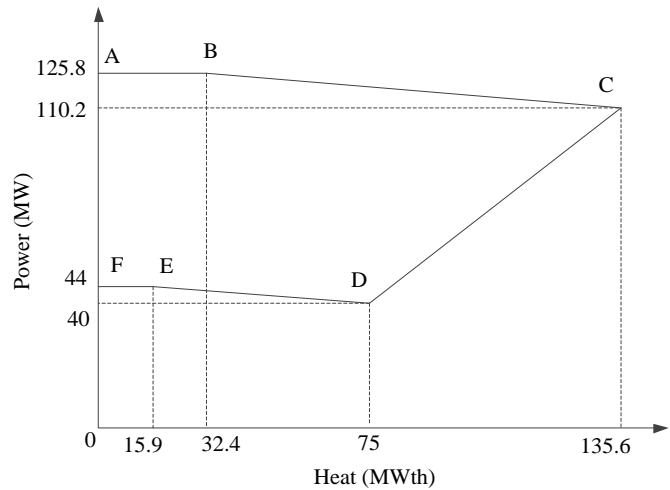

**Fig. 24. FOR of cogeneration unit 2 in test system 2.**

**References**

[1] Sadeghian H R, Ardehali M M. A novel approach for optimal economic dispatch scheduling of integrated combined heat and power systems for maximum economic profit and minimum environmental emissions based on Benders decomposition. Energy 2016, 102:10-23.

[2] Karlsson J, Brunzell L, Venkatesh G. Material-flow analysis, energy analysis, and partial environmental-LCA of a district-heating combined heat and power plant in Sweden. Energy 2018, 144: 31-40.

[3] Zidan A, Gabbar H A, Eldessouky A. Optimal planning of combined heat and power systems within microgrids. Energy 2015, 93: 235-244.

[4] Basu M. Group search optimization for combined heat and power economic dispatch. Int J Electr Power Energy Syst 2016, 78: 138-147.

[5] Rong A, Figueira J R, Lahdelma R. An efficient algorithm for bi-objective combined heat and power production planning under the emission trading scheme. Energy Convers Manage 2014, 88(11):525-534.

[6] Haakana J, Tikka V, Lassila J, et al. Methodology to analyze combined heat and power plant operation considering electricity reserve market opportunities. Energy 2017, 127:408-418.

[7] Montecucco A, Siviter J, Knox A R. Combined heat and power system for stoves with thermoelectric generators. Appl Energy 2017; 185:1336-1342.

[8] Murugan S, Horák B. A review of micro combined heat and power systems for residential applications. Renew Sust Energy Rev 2016; 64: 144-162.

[9] Wang L, Singh C. Stochastic combined heat and power dispatch based on multi-objective particle swarm optimization. Int J Electr Power Energy Syst 2008; 30(3): 226-234.

[10] Shang C, Srinivasan D, Reindl T. Generation and storage scheduling of combined heat and power. Energy 2017; 124: 693-705.

[11] Secui DC. Large-scale multi-area economic/emission dispatch based on a new symbiotic organisms search algorithm. Energy Convers Manage 2017; 154: 203-223.

[12] Niknam T, Azizipanah-Abarghooee R, Roosta A, et al. A new multi-objective reserve constrained combined heat and power dynamic economic emission dispatch. Energy 2012; 42(1): 530-545.

[13] Basu M. Combined heat and power economic emission dispatch using nondominated sorting genetic algorithm-II. Int J Electr Power Energy Syst 2013; 53(1):135-141.

[14] Shi B, Yan L X, Wu W. Multi-objective optimization for combined heat and power economic dispatch with power transmission loss and emission reduction. Energy 2013; 56: 135-143.

[15] Shaabani YA, Seifi AR, Kouhanjani MJ. Stochastic multi-objective optimization of combined heat and power economic/emission dispatch. Energy 2017; 141: 1892-1904.

[16] Ahmadi A, Moghimi H, Nezhad AE, et al. Multi-objective economic emission dispatch considering combined heat and power by normal boundary intersection method. Electr Power Syst Res 2015; 129: 32-43.

[17] Nazari-Heris M, Mohammadi-Ivatloo B, Gharehpetian GB. A comprehensive review of heuristic optimization algorithms for optimal combined heat and power dispatch from economic and environmental perspectives. Renew Sust Energy Rev 2017; 81: 2128-2143.

[18] Li Y, Li Y, Li G, et al. Two-stage multi-objective OPF for AC/DC grids with VSC-HVDC: Incorporating decisions analysis into optimization process. Energy 2018; 147: 286-296.

[19] Walters DC, Sheble GB. Genetic algorithm solution of economic dispatch with valve-point loadings. IEEE Trans Power Syst 1993; 8 (3): 1325-1331.

[20] Niknam T, Mojarrad HD, Meymand HZ. A novel hybrid particle swarm optimization for economic dispatch with valve-point loading effects. Energy Convers Manage 2011; 52(4): 1800-1809.

[21] Abdollahi E, Wang H, Lahdelma R. An optimization method for multi-area combined heat and power production with power transmission network. Appl Energy 2016; 168: 248-256.

[22] Khan NA, Awan AB, Mahmood A, et al. Combined emission economic dispatch of power system including solar photo voltaic generation. Energy Convers Manage 2015; 92:82-91.

[23] Yuan Y, Xu H, Wang B, et al. A new dominance relation-based evolutionary algorithm for many-objective optimization. IEEE Trans Evol Comput 2016; 20(1): 16-37.

[24] Tahboub AM, Pandi VR, Zeineldin HH. Distribution system reconfiguration for annual energy loss reduction considering variable distributed generation profiles. IEEE Trans Power Del 2015; 30(4): 1677-1685.

[25] Mahela OP, Shaik AG. Power quality recognition in distribution system with solar energy penetration using S-transform and Fuzzy C-means clustering. Renew Energ 2017; 106:37-51.

[26] Li H, You S, Zhang H, et al. Modelling of AQI related to building space heating energy demand based on big data analytics. Appl Energy 2017; 203: 57-71.

[27] Li MQ, Yang SX, Liu XH. Diversity comparison of Pareto front approximations in many-objective optimization. IEEE Trans Cybern, 2014; 44(12): 2568-2584.

[28] Deb K, Jain H. An evolutionary many-objective optimization algorithm using reference-point-based nondominated sorting approach, part I: Solving problems with box constraints. IEEE Trans Evol Comput 2014; 18(4): 577-601.

[29] Lu H, Zhang MM, Fei ZM, and Mao KF. Multi-objective energy consumption scheduling in smart grid based on Tchebycheff decomposition. IEEE Trans Smart Grid 2015; 6(6): 2869-2883.

[30] Deb K, Pratap A, Agarwal S, et al. A fast and elitist multiobjective genetic algorithm: NSGA-II. IEEE Trans Evol Comput 2002; 6(2): 182-197.

[31] Vasebi A, Fesanghary M, Bathaee SMT. Combined heat and power economic dispatch by harmony search algorithm. Int J Electr Power Energy Syst 2007; 29(10): 713-719.

[32] Abdolmohammadi HR, Kazemi A. A Benders decomposition approach for a combined heat and power economic dispatch. Energy Convers Manage 2013; 71(1): 21-31.

[33] Ahmadi A, Ahmadi MR. Comment on "Multi-objective optimization for combined heat and power economic dispatch with power transmission loss and emission reduction" Shi B, Yan LX, Wu W [Energy 2013; 56: 226–34]. Energy 2014; 64: 1-2.

[34] Basu M. Artificial immune system for combined heat and power economic dispatch. Int J Electr Power Energy Syst 2012; 43(1):1-5.

[35] Basu M. Bee colony optimization for combined heat and power economic dispatch. Expert Syst Appl 2011; 38(11):13527-13531.

[36] Gu X, Li Y, Jia J. Feature selection for transient stability assessment based on kernelized fuzzy rough sets and memetic algorithm. Int J Electr Power Energy Syst 2015; 64: 664-670.

[37] Haghrah A, Nazari-Heris M, Mohammadi-ivatloo B. Solving combined heat and power economic dispatch problem using real coded genetic algorithm with improved Mühlenbein mutation. Appl Therm Eng 2016; 99: 465–75.

[38] Mohammadi-Ivatloo B, Moradi-Dalvand M, Rabiee A. Combined heat and power economic dispatch problem solution using particle swarm optimization with time varying acceleration coefficients. Electr Power Syst Res 2013; 95:9–18.

[39] Wong KP, Algie C. Evolutionary programming approach for combined heat and power dispatch. Electr Power Syst Res 2002; 61:227–32.

[40] Jordehi AR. A chaotic artificial immune system optimisation algorithm for solving global continuous optimisation problems. Neural Comput Appl 2015; 26:827-33.

[41] Basu M. Combined heat and power economic dispatch by using differential evolution. Electr Power Comp Syst 2010; 38: 996-1004.

[42] Basu M. Dynamic economic emission dispatch using nondominated sorting genetic algorithm-II. Int J Electr Power Energy Syst 2008; 30(2): 140-149.

[43] Li Y, Yang Z, Li G, Zhao D, Tian W. Optimal scheduling of an isolated microgrid with battery storage considering load and renewable generation uncertainties. IEEE Trans Ind Electron 2018, in press. DOI:10.1109/TIE.2018.2840498.

[44] Li Y, Feng B, Li G, Qi J, Zhao D, Mu Y. Optimal distributed generation planning in active distribution networks considering integration of energy storage. Appl Energy 2018; 210, 1073-1081.